\documentclass[11pt]{article}

\usepackage{amssymb}
\usepackage{amsmath}
\usepackage{amsthm}
\usepackage{graphicx}
\usepackage{accents}
\usepackage[margin=0.65cm]{caption}
\usepackage[caption=false]{subfig}
\usepackage{enumerate}
\usepackage{hyperref}
\usepackage{siunitx}

\usepackage{tikz}
\usetikzlibrary{fit,calc}
\usetikzlibrary{arrows,shapes,positioning,shadows,trees}
\usetikzlibrary{decorations.pathreplacing}
\usetikzlibrary{arrows.meta}

\newcommand{\R}{\mathbb{R}}
\newcommand{\N}{\mathbb{N}}
\newcommand{\PO}{\mathbb{P}}
\newcommand{\EE}{\mathbb{E}}
\newcommand{\cB}{{\cal{B}}}
\newcommand{\cC}{{\cal{C}}}
\newcommand{\cI}{{\cal{I}}}
\newcommand{\cM}{{\cal{M}}}
\newcommand{\cP}{{\cal{P}}}
\newcommand{\cJ}{{\cal{J}}}
\newcommand{\cU}{{\cal{U}}}
\newcommand{\hd}{{\hat{d}}}
\newcommand{\bfi}{{\bf i}}

\author{J.~Lang, P.~Domschke, E.~Strauch}
\title{Adaptive Single- and Multilevel Stochastic Collocation Methods
for Uncertain Gas Transport in Large-Scale Networks}
\author{
Jens Lang\footnote{Corresponding author, ORCID: 0000-0003-4603-6554}\\
{\small \it Technische Universit\"at Darmstadt}\\
{\small \it Dolivostra{\ss}e 15, 64293 Darmstadt, Germany}\\[2mm]
Pia Domschke\\
{\small \it Frankfurt School of Finance \& Management}\\
{\small \it Adickesallee 32-34, 60322 Frankfurt am Main, Germany}\\[2mm]
Elisa Strauch\\
{\small \it Technische Universit\"at Darmstadt} \\
{\small \it Dolivostra{\ss}e 15, 64293 Darmstadt, Germany}\\
}
\date{December 7, 2020}
\begin{document}
\maketitle

\abstract{In this paper, we are concerned with the quantification of
uncertainties that arise from intra-day oscillations in the
demand for natural gas transported through large-scale networks. The
short-term transient dynamics of the gas flow is modelled by a hierarchy
of hyperbolic systems of balance laws based on the isentropic Euler equations.
We extend a novel adaptive strategy for solving elliptic PDEs with random data,
recently proposed and analysed by Lang, Scheichl, and Silvester
[J. Comput. Phys., 419:109692, 2020], to uncertain gas transport problems.
Sample-dependent adaptive meshes and a model refinement in the physical space
is combined with adaptive anisotropic sparse Smolyak grids in the stochastic
space. A single-level approach which balances the discretization errors of
the physical and stochastic approximations and
a multilevel approach which additionally minimizes the computational costs
are considered. Two examples taken from a public gas library demonstrate the
reliability of the error control of expectations calculated from
random quantities of interest, and the further use of
stochastic interpolants to, e.g., approximate probability density functions
of minimum and maximum pressure values at the exits of the network.
}

\newpage

\section{Introduction}
\label{lds_sec:intro}
The role of natural gas transport through large-scale networks has been rapidly increased
through the ongoing replacement of traditional energy production by coal fired and
nuclear plants with gas consuming facilities. The safekeeping of energy security and the development of clean energy to meet environmental demands have generated a significant
increase in gas consumption for electric power stations in the last decade. The future
energy mix will mainly be based on low-carbon and regenerative energy and natural gas is considered as a bridging combustible resource to achieve this goal. The seasonally
fluctuating disposability of wind and solar resources causes a growing variability in
electricity production and hence also in the demands of gas transportation by pipelines. The resulting intra-day uncertain oscillations in demand for natural gas leads to new challenges
for computer based modelling and control of gas pipeline operations. Here, an increasing focus
lies on the short-term transient dynamics of gas flow. Operators have to responsively control varying loads to realize a reliable operational management for both gas and electricity delivery systems. These challenging new conditions demand advanced decision tools based on reliable transient simulation and uncertainty quantification taking into account serious operating restrictions.

In this paper, we propose a novel computational approach for the reliable quantification of the transport of uncertainties through a network of gas pipelines. It extends an adaptive multilevel stochastic collocation method recently developed in \cite{LangScheichlSilvester2020} for elliptic partial differential equations with random data to systems of hyperbolic balance laws with
uncertain initial and boundary conditions. We have been developing in-house software tools
for fast and reliable transient simulation and continuous optimization of large-scale gas networks
over the last decade \cite{DomschkeDuaStolwijkEtAl2018,DomschkeKolbLang2010,
DomschkeKolbLang2011b,DomschkeKolbLang2011a,DomschkeKolbLang2015}.
Exemplarily, here we will investigate the important task of safely driving a stationary running system into a newly desired system defined by uncertain gas nominations at delivery points of the network. To be usable in a real-time application of risk analysis and reliability assessment of gas delivery, we have designed our method to meet user-defined accuracies while keeping the computing time for large-scale gas networks at a moderate level. It offers also the opportunity to be integrated in a probabilistic constrained optimization approach \cite{SchusterStrauchGugatLang2020}.

We will consider the following one-dimensional parameterized hyperbolic system of balance laws
on a set of gas pipes $\Omega_j$, $j=1,\ldots,M$, with random initial and boundary data:
\begin{eqnarray}
\partial_t U^{(j)}(x,t,y) + \partial_x F_{m_j}(U^{(j)}(x,t,y)) &=& G_{m_j}(x,t,U^{(j)}(x,t,y)),
\label{lds_eqs:models}\\
U^{(j)}(x,y,0) &=& U_0^{(j)}(x,y),\label{lds_eqs:initial}\\
B(U^{(j)}(x_b,y)) &=& H(x_b,t,y),\quad b\in\cB,\label{lds_eqs:boundary}\\[1mm]
\Phi(U^{(1)}(x_i,y,t),\ldots,U^{(N)}(x_i,y,t))&=&\Pi(x_i,t),\quad i\in\cC,\label{lds_eqs:coupling}
\end{eqnarray}
where the solutions are represented as $U^{(j)}(x,t,y):D^{(j)}\times\Gamma\rightarrow \R^2$
with the deterministic physical domain $D^{(j)}:=\Omega_j\times\R^+$ and
$\Gamma=\Gamma_1\times\Gamma_2\times \cdots \Gamma_N$ being a stochastic parameter space of
finite dimension $N$ (finite noise assumption). The component parameters $y_1,\ldots,y_N$ will
be associated with independent random variables that have a joint probability density
function $\hd(y)=\Pi_{n=1}^{N}\hd_n(y_n)\in L^\infty(\Gamma)$ such that $\hd_n:[-1,1]\rightarrow\R$.
Typically, gas pipeline systems are buried underground and hence temperature differences between a pipe segment and the ground can be neglected in practice. It is therefore standard to consider an isothermal process without a conservation law for the energy, i.e., $U^{(j)}$ is the vector of density $\rho$ and momentum $\rho v$ for each pipe with $v$ being the velocity.

The index sets $\cB$ and $\cC$ in (\ref{lds_eqs:boundary}) and (\ref{lds_eqs:coupling})
describe the indices of the boundary and the coupling nodes, respectively. Boundaries in
gas networks are sources, where gas is injected into the pipeline system, and exits, where
it is taken out by consumers. The modelling of connected pipes, flow at junctions, and the
pressure increase caused by a compressor leads to certain coupling conditions in (\ref{lds_eqs:coupling})
at inner nodes. We ensure conservation of mass and claim the equality of pressure, except for compressors, where the time-dependent term $\Pi(\cdot,t)$ represents the pressure jump that is realised by the compression process. The pressure is calculated from the equation of state for real gases, $p=\rho z(p)RT$,
with compressibility factor $z(p)\in(0,1)$.

We also allow for different gas transport models in each pipe. They are identified by
the parameters $m_j\in\cM:=\{\cM_1,\cM_2,\cM_3\}$ in (\ref{lds_eqs:models}) representing a
whole hierarchy of models with
decreasing fidelity. In our applications, we use the nonlinear isothermal Euler equations as
$\cM_1$, its semilinear approximation as $\cM_2$, and a quasi-stationary model as $\cM_3$. They
will be described in more detail later on.

Let $U=(U^{(1)},\ldots,U^{(M)})$ and $X=C([0,T];L^1(\Omega_1))\times\cdots\times C([0,T];L^1(\Omega_M))$. Throughout this paper, we assume that there is a unique weak entropy solution
$U(\cdot,\cdot,y)\in X$ of the gas flow problem (\ref{lds_eqs:models})-(\ref{lds_eqs:coupling}) for all $y\in\Gamma$. For uncertainty quantification in gas network applications,
it is more natural to consider a functional $\psi(U)$ of the solution $U$ instead of the solution itself. Thus, suppose a possibly nonlinear functional (or quantity of interest) $\psi:X\rightarrow\R$ with $\psi(0)=0$ is given. The
standard collocation method is based on a set of deterministic sample points $\{y^{(q)}\}_{q=1,\ldots,Q}$
in $\Gamma$, chosen to compute independent, finite dimensional space-time approximations $U_h(y^{(q)})\approx U(y^{(q)})$. These approximations are used to construct a single-level interpolant
\begin{equation}
\label{lds_qoi}
\Psi_{Q,h}^{(SL)}(y) = \cI_Q[\psi(U_h)](y) = \sum_{q=1}^Q \psi_q\phi_q(y)
\end{equation}
for the function $\psi(U)$ in the polynomial space with basis functions $\phi_q$,
\begin{equation}
\cP_Q=\text{span}\{\phi_q\}_{q=1,\ldots,Q}\subset L^2_\hd(\Gamma):=
\{\phi:\Gamma\rightarrow\R\text{ s.t. }\int_\Gamma\phi^2(y)\hd(y)dy < \infty\}\,.
\end{equation}
The
interpolation conditions $\cI[\psi(U_h)](y^{(q)})=\Psi(U_h(y^{(q)}))$ for $q=1,\ldots,Q$ determine the
coefficients $\psi_q$. The quality of the interpolation process depends on the accuracy of the space-time
approximations $U_h(y^{(q)})$, the regularity of the solution with respect to the stochastic parameters $y$, and on the number of collocation points $Q$, which grows rapidly with increasing stochastic dimension $N$.
The interpolant $\Psi_{Q,h}^{(SL)}(y)$, also called \textit{response surface approximation}, can be
used to directly calculate moments such as expectation and variance. Since its evaluation is extremely cheap, it also forms the basis for approximating its probability density function by a kernel density
estimator and determining the practically relevant probability that $\Psi_{Q,h}^{(SL)}(y)$ lies in a certain
interval over the whole time horizon. We will apply this approach to check the validity that gas is
delivered in a pressure range stipulated in a contract between gas company and consumer.

The uncertainties in the initial and boundary data in (\ref{lds_eqs:initial}) and (\ref{lds_eqs:boundary})
result in a propagation of uncertainties in the functional $\psi(U)$. It is essential in nowadays natural
gas transport through large networks that operators apply a reliable operational management to guarantee
a \textit{sufficiently smooth} gas flow, respecting at the same time operating limits of compressors and
pressure constraints inside the pipes in a safe manner. There is always a safety factor that prevents the
whole transport system to really hit the limits. Therefore, we may assume appropriate regularity of the solution in the random space in order to ensure a fast convergence of the global approximation polynomials
$\phi_q(y)$ in (\ref{lds_qoi}). Exemplarily, we will investigate the influence of uncertain gas demand when safely driving a stationary running system into a newly desired system defined by shifted gas nominations at the delivery points of the network.

There are two main alternative approaches to stochastic collocation: Monte Carlo sampling and stochastic
Galerkin methods. A detailed discussion of comparative advantages and disadvantages in the context of
hyperbolic systems of conservation laws is given in \cite{AbgrallMishra2017}, see also \cite{BijlLucorMishraSchwab2014,GhanemHigdoOwhadi2016} for a general overview on uncertainty quantification in solutions of more general partial differential equations. Monte Carlo methods and
its variants are the most commonly used sampling methods. They are non-intrusive and robust with respect to
lack of regularity, have a dimension-independent convergence rate and offer a trivial parallelization. However, they are not able to exploit any smoothness or special structure in the parameter dependence and their convergence rate is rather low even when Multilevel Monte Carlo methods are applied. Combined
with finite volume discretizations for the physical space, such methods are extensively investigated
in \cite{MishraRisebroSchwabTokareva2016,MishraSchwab2012,MishraSchwabSukys2012}. Stochastic Galerkin
methods based on generalized polynomial chaos are intrusive and request the solution of heavily
extended systems of conservation laws \cite{PoetteDespresLucor2009,TryonLeMaitreNdjingaErn2010}. Although sparse grids and efficient solvers for block-structured linear systems are used, the computational costs in general are formidable. Recently, an intrusive polynomial moment method which is competitive with
non-intrusive collocation methods has been proposed in \cite{KuschWoltersFrank2020}. In the presence of discontinuities in the random space,
promising semi-intrusive approaches are provided by the stochastic finite volume method
\cite{AbgrallCongedo2013} and a novel hierarchical basis weighted essentially non-oscillatory interpolation method
\cite{Kolb2018}.

The paper is organised as follows. In section 2, we describe the single-level approach and
especially focus on the main ingredients for the adaptive solvers in the physical and
parameter space. The extension to the multilevel approach is explained in section 3, where
we also give asymptotic rates for the complexity of the algorithm. In section 4, two
examples based on networks from a public gas library are investigated to demonstrate the
efficiency and potential of the fully adaptive collocation method. We conclude with a summary
and outlook in section 5.

\section{Adaptive Single-Level Approach}
The main advantage of sampling methods is the reuse of an efficient solver for the transient
gas flow through a network in the range of parameters defined by the stochastic space
$\Gamma$. Since the gas transport through a complex network may be very dynamic and thus changes
both in space and time, an automatic control of the accuracy of the simulation is mandatory.
In order to further reduce computational costs, adjusting the transport model in each pipe
according to the time-dependent dynamics has proven to be very attractive. As a rule of thumb,
the most complex nonlinear Euler equations ($\cM_1$) should be used when needed and the simplest
algebraic model ($\cM_3$) should be taken whenever possible without loosing too much accuracy.
In a series of papers, we have developed a posteriori error estimates and an overall control
strategy to reduce model and discretization errors up to a user-given tolerance
\cite{DomschkeDuaStolwijkEtAl2018,DomschkeKolbLang2010,
DomschkeKolbLang2011b,DomschkeKolbLang2011a,DomschkeKolbLang2015}. A brief introduction will be given next.

Let a parameter $y\in\Gamma$ be fixed and an initial distribution of gas transport models
$\{m_1,\ldots,m_M\}$ be given. Then, we solve the gas network equations (\ref{lds_eqs:models})-(\ref{lds_eqs:coupling}) by means of an adaptive implicit finite
volume discretization \cite{KolbLangBales2010} applied for each pipe until the estimate of
the error in the functional $\psi(U_h(y))$ is less than a prescribed tolerance $\eta_h>0$.
Here, $h$ refers to resolution in space, time and model hierarchy. To raise efficiency, the
simulation time $[0,T]$ is divided into subintervals $[t_i,t_{i+1}]$, $i=0,\ldots,N_t-1$, of
the same size. We then successively process the classical adaption loop
\begin{equation}
\label{lds_adapt_loop}
\text{SOLVE} \Rightarrow \text{ESTIMATE} \Rightarrow \text{MARK} \Rightarrow\text{REFINE}
\end{equation}
for each of the subintervals such that eventually
\begin{equation}
\label{lds_qoi_acc}
\left| \psi(U(y)) - \psi(U_h(y)) \right| \le c_h(y)\, \left| \sum_{j=1}^M
\left( \eta_{x,j}+\eta_{t,j}+\eta_{m,j} \right) \right| < c_h(y)\cdot\eta_h
\end{equation}
in the second step with a sample-dependent constant $c_h(y)$ that is usually
close to one. The a posteriori error estimators $\eta_{x,j}$, $\eta_{t,j}$,
and $\eta_{m,j}$ for the $j$-th pipe determine the error distribution along the network for the
spatial, temporal and model discretizations. They measure the influence of the
model and the discretisation on the output functional $\psi$ and can be
calculated by using the solutions of adjoint equations. A detailed description
which would go beyond the scope of our paper is given in
\cite[Sect. 2.2]{DomschkeDuaStolwijkEtAl2018}, see also
\cite{DomschkeKolbLang2011b,DomschkeKolbLang2011a}. Polynomial reconstructions in
space and time of appropriate orders are used to compute $\eta_{x,j}$ and $\eta_{t,j}$,
respectively. The model error estimator $\eta_{m,j}$ is derived from the
product of differential terms, representing the difference between models, and the
sensitivities calculated from the adjoint equations.

In our calculations, we use the following model hierarchy:
\begin{itemize}\leftskip15pt
\item $\cM_1$: Nonlinear isothermal Euler equations
\begin{eqnarray*}
\partial_t \rho + \partial_x (\rho v)&=& 0,\\
\partial_t (\rho v) + \partial_x (p+\rho v^2)&=& g(\rho,\rho v),
\end{eqnarray*}
\item $\cM_2$: Semilinear isothermal Euler equations
\begin{eqnarray*}
\partial_t \rho + \partial_x (\rho v)&=& 0,\\
\partial_t (\rho v) + \partial_x p &=& g(\rho,\rho v),
\end{eqnarray*}
\item $\cM_3$: Algebraic isothermal Euler equations
\begin{eqnarray*}
\partial_x (\rho v)&=& 0\\
\partial_x p &=& g(\rho,\rho v)
\end{eqnarray*}
\end{itemize}
with the joint source term $g(\rho,\rho v)=-\lambda\rho v|v|/(2D)$,
where $D$ is the pipe diameter and $\lambda$ the Darcy friction
coefficient. We note that the algebraic model can be analytically solved
in the variables $\rho v$ and $p$. The models are connected at inner nodes,
where we ensure conservation of mass and the equality of pressure. The
latter one is often used in engineering software, but can be also replaced
by the equality of total enthalpy. The interested reader is referred to the
discussion in \cite{MindtLangDomschke2019}. Pipes can also be connected
by valves and compressors. Valves are used to regulate the flow in gas
networks. An open valve is modelled as inner node, whereas $q=0$ is
required at both sides of a closed valve. Compressors compensate for the
pressure loss due to friction in the pipes. The power of
a compressor $c\in\cJ$ that is needed for the compression process is
given by
\begin{equation}
\label{lds_gc}
G_c(U(t)) = c_F\,q_{in}(t)\,z(p_{in}(t))\left(\left(
\frac{p_{out}(t)}{p_{in}(t)}\right)^{\frac{\gamma-1}{\gamma}} - 1 \right)
\end{equation}
with in- and outgoing pressure $p_{in}$, $p_{out}$, and ingoing flow rate
$q_{in}$ \cite{Menon2005}.
The parameter $c_F$ is a compressor specific constant, $\gamma$ the isentropic
coefficient of the gas, and $z\in (0,1)$ the compressibility factor from the
equation of state for real gases. In our application, we use the specific energy
consumption needed by the electric motors to realize all desired compressions
as quantity of interest that drives the adaptation process. It can be estimated by
a quadratic polynomial in $G_c$, i.e., we set
\begin{equation}
\label{lds_psi_uy}
\psi(U(y)) = \alpha\,\sum_{c\in\cJ} \int_{0}^{T} g_{c,0} + g_{c,1}\,G_c(U(y)) +
g_{c,2}\,G_c^2(U(y))\,dt
\end{equation}
with given compressor-dependent constants $g_{c,i}\in\R$ and a scaling
factor $\alpha>0$.

The complex task in the step
MARK (for refinement) of finding an optimal refinement strategy that
combines the three types of adaptivity is a generalisation of the
unbounded knapsack problem, which is NP-hard. A good approximation
can be found by a greedy-like refinement strategy as investigated in
\cite{DomschkeDuaStolwijkEtAl2018}. It leads to considerable computational
savings without compromising on the simulation accuracy. Eventually,
we have an adaptive black box solver - our working horse - at hand
that, once a random parameter $y\in\Gamma$ and a specific tolerance $\eta_h$
have been chosen, delivers a numerical approximation $U_h(y)$ such that
the accuracy requirement (\ref{lds_qoi_acc}) is satisfied for
\begin{equation}
\label{lds_adapt_network}
\psi(U_h(y)) = \text{ANet}(y,\eta_h)\,.
\end{equation}
Working close to the asymptotic regime, we can assume that the adaptive
algorithm converges for fixed $y\in\Gamma$ and $\eta_h\rightarrow 0$.

Starting from the pointwise error estimate (\ref{lds_qoi_acc}) and supposing
bounded first moments of $c_h(y)$, we directly get the following error bound:
\begin{equation}
\label{lds_err_h}
\left| \EE [\psi(U(y))-\psi(U_h(y))] \right| :=
\left| \int_{\Gamma} \left( \psi(U(y))-\psi(U_h(y)) \right) \hd(y)\,dy \right|
 \le C_h\cdot \eta_h
\end{equation}
with a constant
\begin{equation}
\label{lds_err_h_const}
C_h := \int_\Gamma c_h(y)\,\hd(y)\,dy
\end{equation}
that does not depend on $y$.

We will now discuss the control of the error for the adaptive
stochastic collocation method. Let us assume
$\psi(U)\in C^0(\Gamma,\R)$ and consider the interpolation operator
$\cI_Q:C^0(\Gamma)\rightarrow L^2_\hd(\Gamma)$ from (\ref{lds_qoi}). This
operator is constructed by a hierarchical sequence of one-dimensional
Lagrange interpolation operators, using the anisotropic Smolyak algorithm
as introduced in \cite{GerstnerGriebel2003}. It reads
\begin{equation}
\label{lds_sparse_grid_intpl}
\begin{array}{rll}
\cI_Q[\psi(U_h)](y) &=& \sum_{\bfi\in I} \triangle^{m(\bfi)}[\psi(U_h)](y)\\[2mm]
&:=&\sum_{\bfi\in I} \bigotimes_{n=1}^{N}
\left( \cI_n^{m(i_n)}[\psi(U_h)](y) - \cI_n^{m(i_n-1)}[\psi(U_h)](y)\right)
\end{array}
\end{equation}
with multi-indices $\bfi=(i_1,\ldots,i_N)\in I\subset\N^N_+$,
$m(\bfi)=(m(i_1),\ldots,m(i_N))$, and
univariate polynomial interpolation operators
$\cI_n^{m(i_n)}:C^0(\Gamma_n)\rightarrow\PO_{m(i_n)-1}$. These operators
use $m(i_n)$ collocation points to construct a polynomial interpolant
in $y_n\in\Gamma_n$ of degree at most $m(i_n)-1$. The operators
$\triangle^{m(\bfi)}$ are often referred to as hierarchical
surplus operators. The function
$m$ has to satisfy $m(0)=0$, $m(1)=1$, and $m(i)<m(i+1)$. We formally set
$\cI^0_n=0$ for all $n=1,\ldots,N$ and
use the nested sequence of univariate Clenshaw--Curtis nodes with
$m(i)=2^{i-1}+1$ if $i>1$. The index $Q$ in (\ref{lds_sparse_grid_intpl}) is
then the number of all explored quadrature points in $\Gamma$
determined by $m(\bfi)$.

The value of the hierarchical surplus operator $\triangle^{m(\bfi)}$
in (\ref{lds_sparse_grid_intpl}) can be interpreted as profit and therefore
used as error indicator for already computed approximations. Applying once
again the classical adaption loop from (\ref{lds_adapt_loop}), the adaptive
anisotropic Smolyak algorithm computes profits in each step, adds the index
of the highest profit to the index set $m(\bfi)$ and explores admissible
neighbouring indices next. The algorithm stops if the absolute value of the
highest profit is less than a prescribed tolerance, say $\eta_s>0$. Obviously,
the method is dimension adaptive. There is a {\sc Matlab} implementation
\textit{Sparse Grid Kit} available, which can be downloaded from the CSQI
website \cite{BeckNobileTamelliniTempone2011}. Its numerical performance is
discussed in the review paper \cite{TamelliniNobileGuignardTeseiSprungk2017}.

Following this adaptive methodology, we get an error estimate
\begin{equation}
\label{lds_err_s}
\left| \EE \left[\psi(U_h(y))-\cI_Q[\psi(U_h)](y)\right]\right| \le C_s\cdot\eta_s
\end{equation}
with a constant $C_s>0$. We assume that $C_s$ does not depend on $h$. If we
now split the overall error into the sum of a physical error resulting
from the chosen resolution in space, time and model hierarchy, and a stochastic
interpolation error, then using the inequalities (\ref{lds_err_h}), (\ref{lds_err_s})
and the triangle inequality yields the final estimate
\begin{equation}
\label{lds_err_sl}
\begin{array}{rll}
&& \hspace{-1cm}\left| \EE \left[ \psi(U(y))-\cI_Q[\psi(U_h)](y) \right] \right|  \\[2mm]
& \le & \left| \EE [\psi(U(y))-\psi(U_h(y))] \right| +
\left| \EE \left[ \psi(U_h(y))-\cI_Q[\psi(U_h)](y)\right]\right| \\[2mm]
& \le & C_h\cdot\eta_h + C_s\cdot\eta_s\,.
\end{array}
\end{equation}
Let $\varepsilon>0$ be a user-prescribed tolerance for the error on the
left-hand side. Then the usual strategy to balance both the physical
and the stochastic error on the right-hand side is to choose the individual
tolerances as $\eta_h=\varepsilon/(2C_h)$ and $\eta_s=\varepsilon/(2C_s)$.
Finally, the adaptive Smolyak algorithm is called with the tolerance $\eta_s$,
where for each chosen sample point $y\in\Gamma$, the black box solver in
(\ref{lds_adapt_network}) runs with ANet($y,\eta_h$), resulting in a
sample-adaptive resolution in the physical space. The algorithm is illustrated
in Tab.~\ref{lds_tab:1}.
\begin{table}[!t]
\centering
\parbox{14cm}{
\caption{\small Algorithm to approximate solution functionals $\psi(U)$ by an
adaptive single-level stochastic collocation method.}
\label{lds_tab:1}}
\begin{tabular}{p{0.5cm}p{10.5cm}}
\hline\noalign{\smallskip}
\multicolumn{2}{l}{Algorithm: Adaptive Single-Level Stochastic Collocation Method} \\
\noalign{\smallskip}\hline\noalign{\smallskip}
1. & Given $\varepsilon$, estimate $C_h$, $C_s$, and set $\eta_h:=\varepsilon/(2C_h)$,
$\eta_s:=\varepsilon/(2C_s)$.\\[1mm]
2. & Compute $\EE[\cI_Q[\psi(U_h)]]$ := ASmol(ANet($y,\eta_h$),$\eta_s$).\\
\noalign{\smallskip}\hline\noalign{\smallskip}
\end{tabular}
\end{table}
\section{Adaptive Multilevel Approach}
Next, we consider an adaptive multilevel approach in order to enhance the
efficiency of the uncertainty quantification further. First multilevel strategies
in the context of Monte Carlo methods were independently proposed as an abstract
variance reduction technique in \cite{Giles2008,Heinrich2001}. Extensions to
uncertainty quantification were developed in \cite{BarthSchwabZollinger2011,CliffeGilesScheichlTeckentrup2011}. Later on, they
also entered the field of stochastic collocation methods
\cite{LangScheichlSilvester2020,TeckentrupJantschWebsterGunzburger2012,ZhuLinebargerXiu2017}.
The methodology in this paper can be viewed as an extension of the adaptive multilevel
stochastic collocation method developed for elliptic PDEs with random data in \cite{LangScheichlSilvester2020} to the hyperbolic case, where a sample-dependent hierarchy of spatial approximations is replaced by a more sophisticated space-time-model hierarchy.

Let a sequence $\{\eta_{h_k}\}_{k=0,\ldots,K}$ of tolerances with
\begin{equation}
1 \ge \eta_{h_0} > \eta_{h_1} > \ldots > \eta_{h_K} > 0
\end{equation}
be given. Each $h_k$ refers to a certain resolution in space, time and model hierarchy
such that for any solution $U_{h_k}(y)$ with $y\in\Gamma$ it holds
\begin{equation}
\label{lds_err_allh}
\left| \EE [\psi(U(y))-\psi(U_{h_k}(y))] \right| \le C_H\cdot \eta_{h_k},\quad k=0,\ldots,K,
\end{equation}
with a constant $C_H := \max_{k=0,\ldots,K} C_{h_k}$ that does not depend on $y$. The constants
$C_{h_k}$ are defined in (\ref{lds_err_h_const}) with $h=h_k$. We consider now a second
family of (stochastic) tolerances $\{\eta_{s_k}\}_{k=0,\ldots,K}$ and assume that there exits
numbers $Q_k$, $k=0,1,\ldots,K$ and a positive constant $C_Y$ not depending
on $k$ such that
\begin{equation}
\label{lds_err_allk}
\left| \EE \left[ \left( \psi(U_{h_k}) - \psi(U_{h_{k-1}}) \right) -
\cI_{Q_{K-k}} \left[ \psi(U_{h_k}) - \psi(U_{h_{k-1}}) \right]\right]\right|
\le C_Y\cdot\eta_{s_{K-k}}
\end{equation}
for $k=0,1,\ldots,K$. Here, we formally set $\psi(U_{h_{-1}}):=0$. Observe that with
increasing index $k$, the differences $|\psi(U_{h_k})(y) - \psi(U_{h_{k-1}})(y)|$
decrease and hence the number of collocation points $Q_{K-k}$ necessary to achieve the tolerance
$\eta_{s_{K-k}}$ gets smaller. Consequently, less samples on fine meshes and with high fidelity
models are needed to achieve the overall tolerance, which is the main motivation for the use
of a multilevel approach.

Using a telescopic sum of single-level interpolants, we construct a multilevel interpolant for
the functional $\psi(U)$ through
\begin{equation}
\begin{array}{rll}
\Psi^{(ML)}_{K}(y)&:=&\sum_{k=0}^K \left( \Psi^{(SL)}_{Q_{K-k},h_k}(y) -
\Psi^{(SL)}_{Q_{K-k},h_{k-1}}(y) \right)\\[4mm]
&=& \sum_{k=0}^K \cI_{Q_{K-k}}\left[\psi(U_{h_k}) - \psi(U_{h_{k-1}}) \right](y).
\end{array}
\end{equation}
Its error can be estimated by
\begin{equation}
\label{lds_err_ml}
\begin{array}{rll}
&& \hspace{-1cm}\left| \EE [\psi(U(y))-\Psi^{(ML)}_K(y)] \right| \\[2mm]
&\le&
\left| \EE [\psi(U(y))-\psi(U_{h_K}(y))] \right| +
\left| \EE [\psi(U_{h_K}(y))-\Psi^{(ML)}_K(y)] \right| \\[2mm]
&\le& C_H\cdot\eta_{h_K} + C_Y\cdot\sum_{k=0}^K \eta_{s_{K-k}},
\end{array}
\end{equation}
where we have used the identity
$\psi(U_{h_K})=\sum_{k=0,\ldots,K}(\psi(U_{h_k})-\psi(U_{h_{k-1}}))$
and the inequalities (\ref{lds_err_allh}) and (\ref{lds_err_allk}). There
are two different ways to balance the errors on the right-hand side:
(1) set $\eta_{s_k}=C_H\cdot\eta_{h_K}/((K+1)C_Y)$ for all $k=0,\ldots,K$,
which yields $2C_H\cdot\eta_{h_K}$ as upper bound in (\ref{lds_err_sl}), and
(2) choose $\eta_{s_k}$ in such a way that the computational cost is minimized.
We will go for the second option and follow the suggestions given in
\cite{LangScheichlSilvester2020}. Let $W_k$ denote the work (computational cost)
that must be invested to solve the gas network equations for a sample point
$y\in\Gamma$ with accuracy $\eta_{h_k}$. Then we make the following assumptions:
\begin{equation}
\begin{array}{crll}
\text{(A1)}\hspace{1cm} & W_k & \le & C_W\cdot\eta_{h_k}^{-s}, \\[2mm]
\text{(A2)}\hspace{1cm} & C_Y\cdot\eta_{s_{K-k}} & = &
C_I(N)M_{K-k}^{-\mu}\eta_{h_{k-1}},
\end{array}
\end{equation}
for all $k=0,\ldots,K$. Here, we fix
$\eta_{h_{-1}}:=|\EE[\cI_{Q_0}[\psi(U_{h_0})]]|$. The constants $C_W>0$,
$C_I(N)>0$ are independent of $k$, $y$, and the rates $s$ and $\mu$ are
strictly positive. Recall that $N$ is the dimension of the stochastic space.

To achieve an accuracy $\varepsilon>0$ for the multilevel interpolant, i.e.,
\begin{equation}
\left| \EE [\psi(U(y))-\Psi^{(ML)}_K(y)] \right| \le \varepsilon,
\end{equation}
at minimal cost $C_\varepsilon^{(ML)}:=\sum_{k=0,\ldots,K}Q_{K-k}(W_k+W_{k-1})$,
the optimal choice of the stochastic tolerances is given in
\cite[Theorem 2.1]{LangScheichlSilvester2020}. They are
\begin{equation}
\eta_{s_{K-k}} = \left( 2C_YG_K(\mu)\right)^{-1}
\left( F_k(s)\right)^{\frac{\mu}{\mu+1}}\eta_{h_{k-1}}\varepsilon,
\quad k=0,\ldots,K,
\end{equation}
where
\begin{equation}
\begin{array}{rll}
F_0(s) &=& \eta_{h_0}^{-s}\eta_{h_{-1}}^{-1},\\[4mm]
F_k(s) &=& \left( \eta_{h_k}^{-s} +
\eta_{h_{k-1}}^{-s} \right)\eta_{h_{k-1}}^{-1},\quad k=1,\ldots,K,\\[4mm]
G_K(\mu) &=& \sum_{k=0}^{K}\left( F_k(s)\right)^{\frac{\mu}{\mu+1}}\eta_{h_{k-1}}.
\end{array}
\end{equation}
Typically, in practical calculations, a decreasing sequence of tolerances
$\eta_{h_k}=q^k\eta_{h_0}$ with a positive reduction factor $q<1$ is used.
In this case, we can estimate the overall multilevel costs using a standard
construction \cite[Theorem 2.2]{LangScheichlSilvester2020},
\begin{equation}
C_\epsilon^{(ML)} \lesssim \left\{
\begin{array}{ll}
\epsilon^{-\frac{1}{\mu}}                                 & \mbox{if } s\mu<1\\[1mm]
\epsilon^{-\frac{1}{\mu}}|\log\epsilon |^{1+\frac{1}{\mu}}& \mbox{if } s\mu=1\\[1mm]
\epsilon^{-s}                                             & \mbox{if } s\mu>1.
\end{array}
\right.
\end{equation}
\begin{table}[!t]
\centering
\parbox{15cm}{
\caption{\small Algorithm to approximate solution functionals $\psi(U)$ by an
adaptive multilevel stochastic collocation method.}
\label{lds_tab:2}}
\begin{tabular}{p{0.1cm}p{13cm}}
\hline\noalign{\smallskip}
\multicolumn{2}{l}{Algorithm: Adaptive Multilevel Stochastic Collocation Method} \\
\noalign{\smallskip}\hline\noalign{\smallskip}
1. & Given $\varepsilon$, $q$, and $K$, estimate $C_H$, $C_Y$, $s$, $\mu$ and set
$\eta_{h_K}:=\varepsilon/(2C_H)$.\\[1mm]
2. & Set $\eta_{h_k}:=q^{k-K}\eta_{h_K}$ for $k=0,\ldots,K-1$.\\[1mm]
3. & Compute $\eta_{h_{-1}}:=\EE[\cI_{Q_0}[\psi(U_{h_0})]]$
:= ASmol(ANet($y,\eta_{h_0}$),$\eta_{h_0}$).\\[1mm]
4. & Set $\eta_{s_{K-k}}:=\left( 2C_YG_K(\mu)\right)^{-1}
\left( F_k(s)\right)^{\frac{\mu}{\mu+1}}\eta_{h_{k-1}}\varepsilon$ for
$k=0,\ldots,K$.\\[1mm]
5. & Compute $E_0$ := ASmol(ANet($y,\eta_{h_0}$),$\eta_{s_K}$).\\[1mm]
6. & Compute $E_k$ := ASmol(ANet($y,\eta_{h_k}$) - ANet($y,\eta_{h_{k-1}}$),$\eta_{s_{K-k}}$)
for $k=1,\ldots,K$.\\[1mm]
7. & Compute $\EE[\Psi_K^{(ML)}]:=\sum_{k=0,\ldots,K}E_k$.\\
\noalign{\smallskip}\hline\noalign{\smallskip}
\end{tabular}
\end{table}

For the convenience of the reader, we summarize the multilevel algorithm in
Table~\ref{lds_tab:2}. Once the parameters are set, the approach is self-adaptive
in nature. Observe that already computed samples at level $k-1$ can be reused to
compute $E_k$ in step~6. In general, sufficient estimates for the constants
$C_H$, $C_Y$ and the rates $\mu$, $s$ can be derived from the study of a few
samples with relatively coarse resolutions. In any case, these samples can be
reused later on.

\section{Practical Examples from a Gas Library}
Examplarily, we will consider two gas network configurations: \emph{GasLib-11} and
\emph{GasLib-40} from the public gas library \url{gaslib.zib.de}
\cite{SchmidtAsmannBurlacuEtAl2017}. They
are parts of real gas networks in Germany. We have implemented the adaptive approach
described above for the deterministic black box solver \textit{ANet($y,\eta_h$)}
from (\ref{lds_adapt_network}) in our in-house software package {\sc Anaconda}.
More details of the implementation can be found in~\cite{Kolb2011}. The adaptive
stochastic collocation method \textit{ASmol($\cdot,\eta_s$)} was realized by means
of the \textit{Sparse Grid Kit}
developed in {\sc Matlab} \cite{TamelliniNobileGuignardTeseiSprungk2017}.
All calculations have been done with {\sc Matlab} version R2020a on a
Intel(R) Xeon(R) Gold 6130 CPU running at 2.1 GHz.

A common daily operation of gas networks is the smooth transformation of a stationary state
$U_A$, which has worked well for the given demands so far, into a new stationary state $U_B$,
which realizes a change in the gas demand over a couple of hours. This scenario is best
treated by appropriate optimization tools which determine the operating range of all compressor
stations and valves in such a way that, e.g., lower and upper bounds of pressures are
satisfied during the whole time-dependent conversion process. In what follows, we will assume
that a feasible, optimized control, i.e., the operation modes for compressors (pressure jump)
and valves (open or closed),
is already known for this so-called nomination change. Then, we will fix these controls and focus
on the influence of uncertainties in the consumers' demands around state $U_B$ on the compressor costs and the
feasibility of the pressure at which the gas is delivered to the consumers. Typically,
corresponding pressure requirements are regulated in contracts.

\subsection{An Example with 11 Pipes}
\begin{figure}[t]
\centering
\includegraphics[scale=.7]{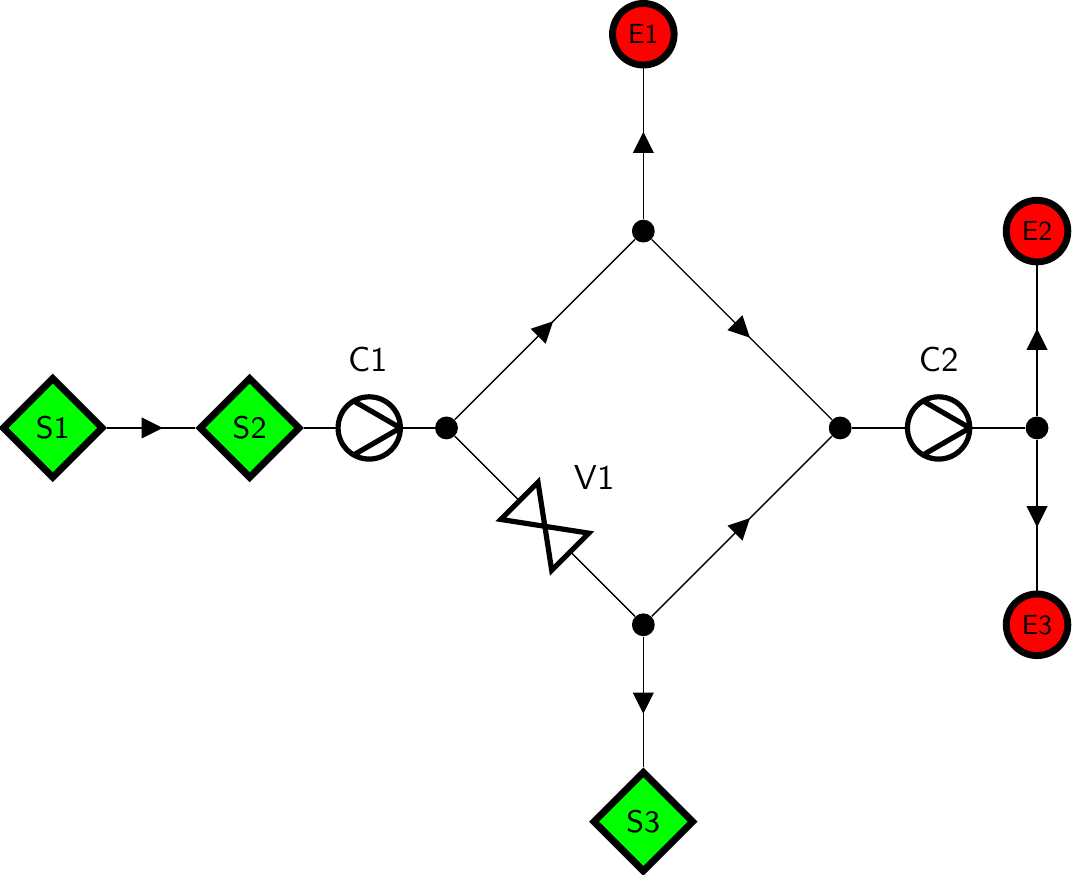}\\[2mm]
\parbox{14cm}{
\caption{\small Schematic description of the network \emph{GasLib-11} with 11 pipes, 2 compressors (C1,C2),
1 valve (V1), 3 sources (green diamonds: S1, S2, S3) and 3 exits (red circles: E1, E2, E3).
The arrows determine the
orientation of the pipes to identify the flow direction by the sign of the velocity.}
\label{lds_fig:gaslib11}
}
\end{figure}

The first example is taken from the \emph{GasLib-11}, which consists of $11$ pipes, $2$ compressors,
$1$ valve, $3$ sources, and $3$ exits, see Fig.~\ref{lds_fig:gaslib11}. The stationary initial state
$U_A$ and the final state $U_B$ are determined by the boundary conditions and controls given in
Tab.~\ref{lds_tab:3}. The simulation is started with $U_0=U_A$. After $4h$, the boundary
values and controls are linearly changed to reach the new conditions defined for $U_B$ at $t=6h$.
The valve is closed
at $t=4.5h$. The simulation time of $24h$ is split into subintervals of $4h$, for which the
classical adaption loop (\ref{lds_adapt_loop}) is processed.
\begin{table}[!h]
\centering
\parbox{14cm}{
\caption{\small\emph{GasLib-11}: Boundary data for sources (S1-S3), exits (E1-E3), and controls for
compressors (C1-C2) and valves (V1) for initial state $U_A$ and final state $U_B$.}
\label{lds_tab:3}}
\begin{tabular}{|p{3.5cm}|ccc||ccc|} \hline
 & \multicolumn{3}{|c||}{State $U_A$} & \multicolumn{3}{|c|}{State $U_B$}\\ \hline\hline
source & S1 & S2 & S3 & S1 & S2 & S3 \\ \hline
pressure [\si{\bar}] & 70.00 & 70.00 & 65.00 & 48.00 & 54.00 & 46.00 \\ \hline\hline
exit &  E1 & E2 & E3 &  E1 & E2 & E3 \\ \hline
volume flow [\si{m^3.s^{-1}}] & 38.22 & 38.22 & 38.22 & 25.48 & 25.48 & 25.48\\ \hline\hline
compressor & C1 & C2 & & C1 & C2 & \\ \hline
pressure jump [\si{\bar}] & 0 & 0 & & 5 & 15 & \\ \hline\hline
valve &  V1 & & &  V1 & &\\ \hline
operation & open &  & & closed & & \\ \hline
\end{tabular}
\end{table}
For the state $U_B$, the volume flows $q_{E}$ at the three exists $E\!=\!E1,E2,E3$,
are \textit{uncertain} due to an individual behaviour
of the consumers and are parameterised by three variables $y=(y_1,y_2,y_3)$, representing the image
of a triple of independent random variables with $y_i\in\cU[-1,1]$. We set
\begin{equation}
q_{Ei}(y_i) = 25.48 + 10\cdot y_i,\quad i=1,2,3.
\end{equation}
According to (\ref{lds_psi_uy}), the quantity of interest $\psi$ is defined by
the specific energy consumption of the compressors,
\begin{equation}
\psi(U(y)) = \alpha\,\sum_{c=C1,C2} \int_{0h}^{24h} g_{c,0} + g_{c,1}\,G_c(U(y)) +
g_{c,2}\,G_c^2(U(y))\,dt
\end{equation}
with $g_{c,0}=5000$, $g_{c,1}=2.5$, $g_{c,2}=0$ for both compressors and
$G_c$ defined in (\ref{lds_gc}). We set the weighting factor $\alpha=10^{-10}$ to bring the
expected value of $\psi(U(y))$ in the order of $0.1$.

\begin{figure}[t]
\centering
\includegraphics[scale=.6]{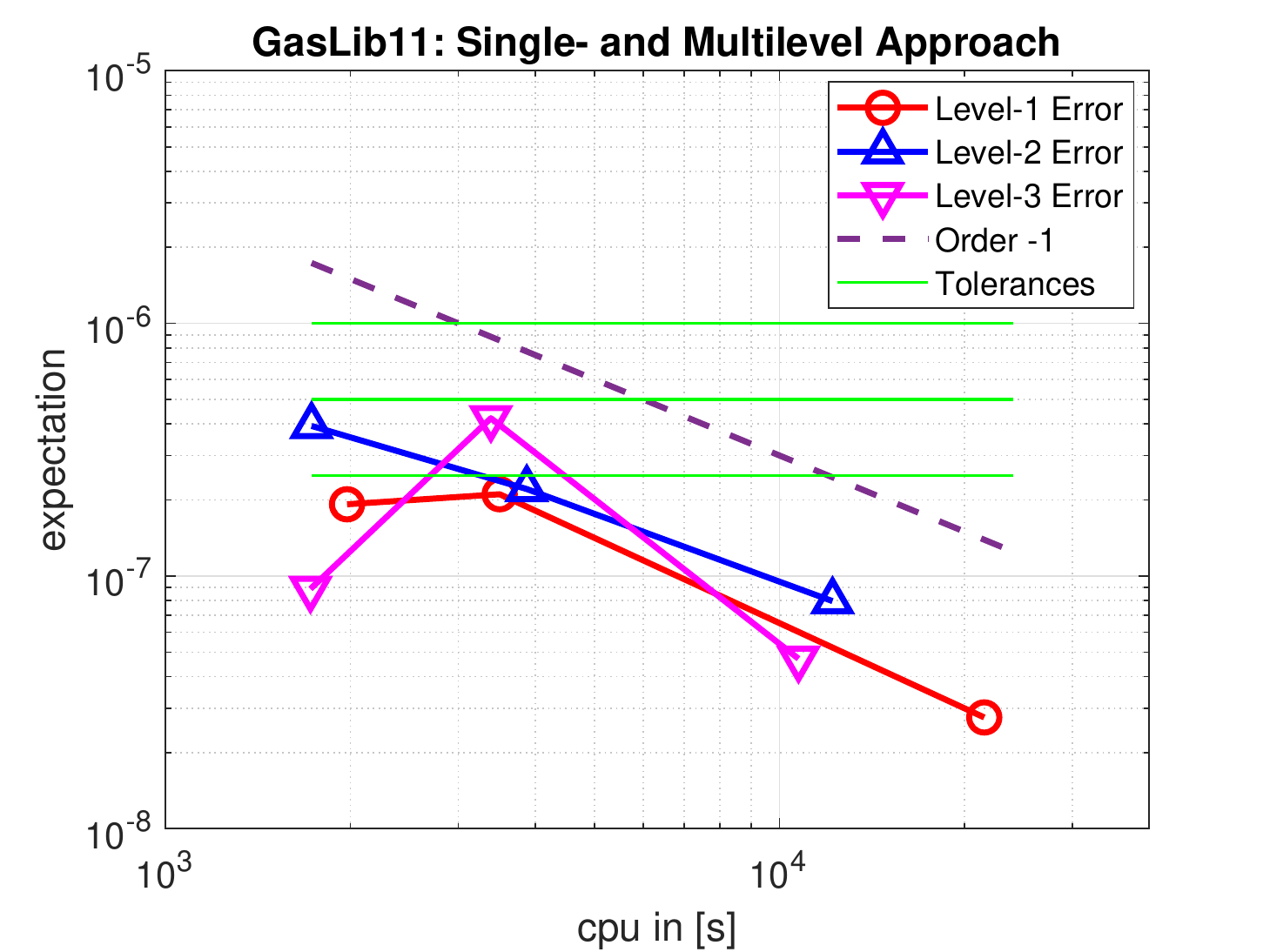}\\[2mm]
\parbox{14cm}{
\caption{\small\emph{GasLib-11}: Errors for the expected values $\EE[\psi_K^{(ML)}]$
and $\EE[\psi^{(SL)}]$ for the three-level (magenta triangles down, $K\!=\!2$),
two-level (blue triangles up, $K\!=\!1$), and one-level (red circles) approach with
adaptive space-time-model discretizations for
$\varepsilon=10^{-6},5\times 10^{-7},2.5\times 10^{-7}$ (green lines). The
accuracy achieved is almost always better than the tolerance. The single-level
and the two-level approach perform quite similar. The three-level approach shows
an irregular behaviour, but also delivers very good results.}
\label{lds_fig:slml11}
}
\end{figure}

In order to start the adaptive stochastic collocation method, we have performed a few
calculations for low tolerances to estimate the parameters in Tab.~\ref{lds_tab:2}.
We found $C_H=C_Y=0.1$, $s=1$, and $\mu=2$ by a least squares fit and appropriate
rounding. Let us now consider
a single-, two, and three-level approach with a reduction factor $q=0.5$. The
overall accuracy requirements are
$\varepsilon=10^{-6},5\times 10^{-7},2.5\times 10^{-7}$, where a reference
solution $\EE[\psi(U)]=0.123765671196008$ is calculated with $\varepsilon=5\times 10^{-8}$.
The results are shown in Fig.~\ref{lds_fig:slml11}. The differences of the methods
are not very high. This is due to the fact that only $25$ collocation points are
sufficient to reach the highest accuracy.

\begin{figure}[!h]
\centering
\includegraphics[scale=.45]{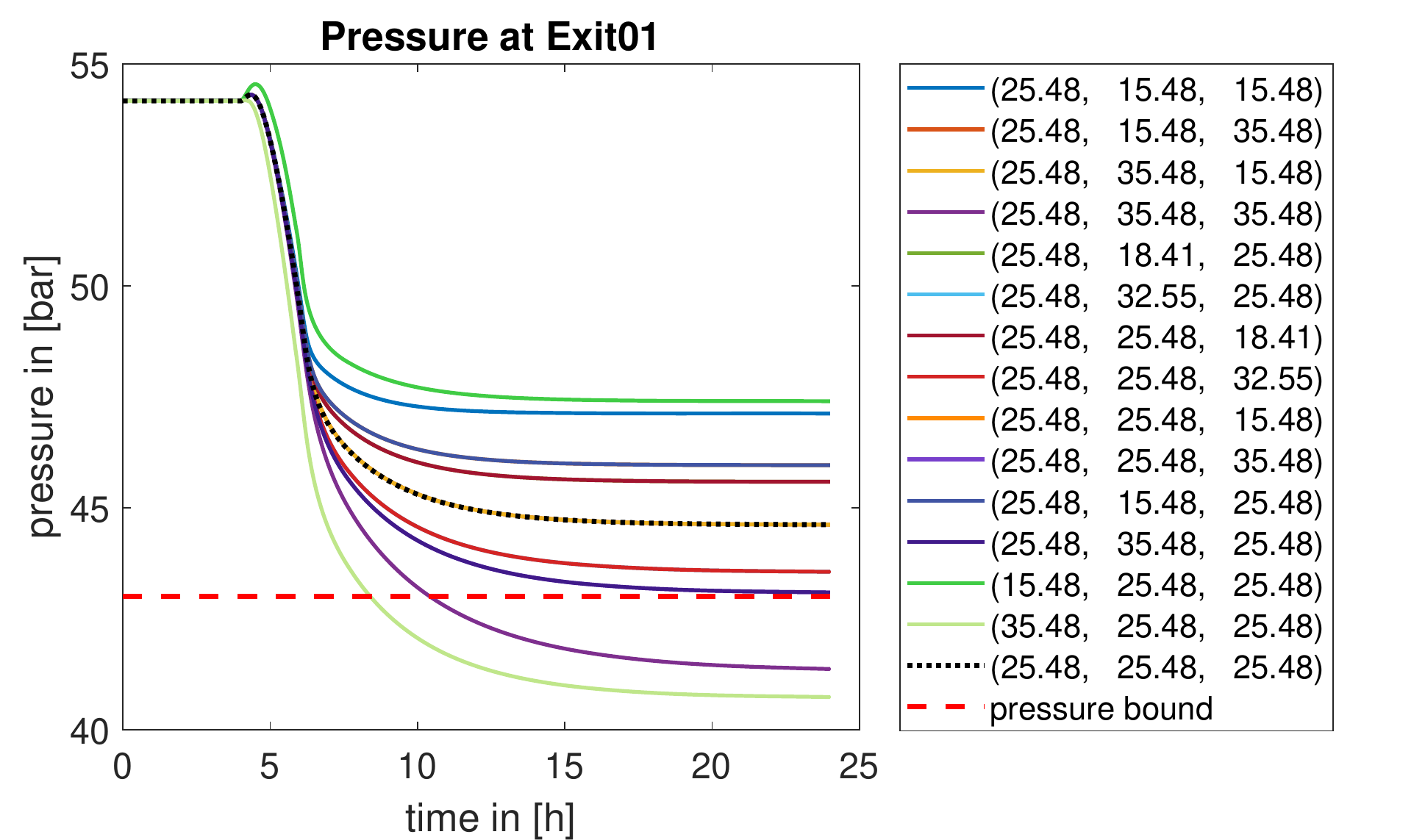}\\[2mm]
\includegraphics[scale=.43]{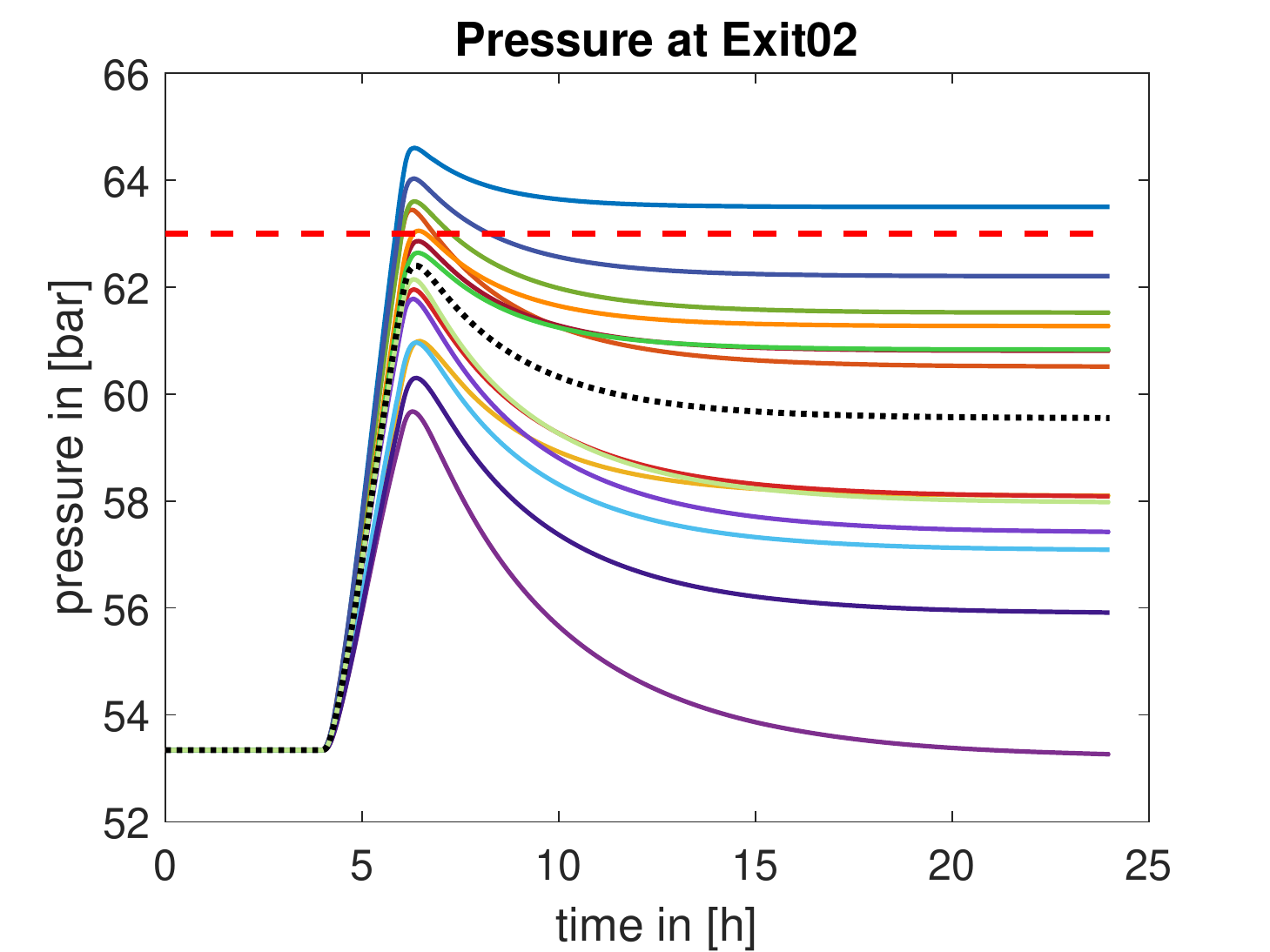}
\includegraphics[scale=.43]{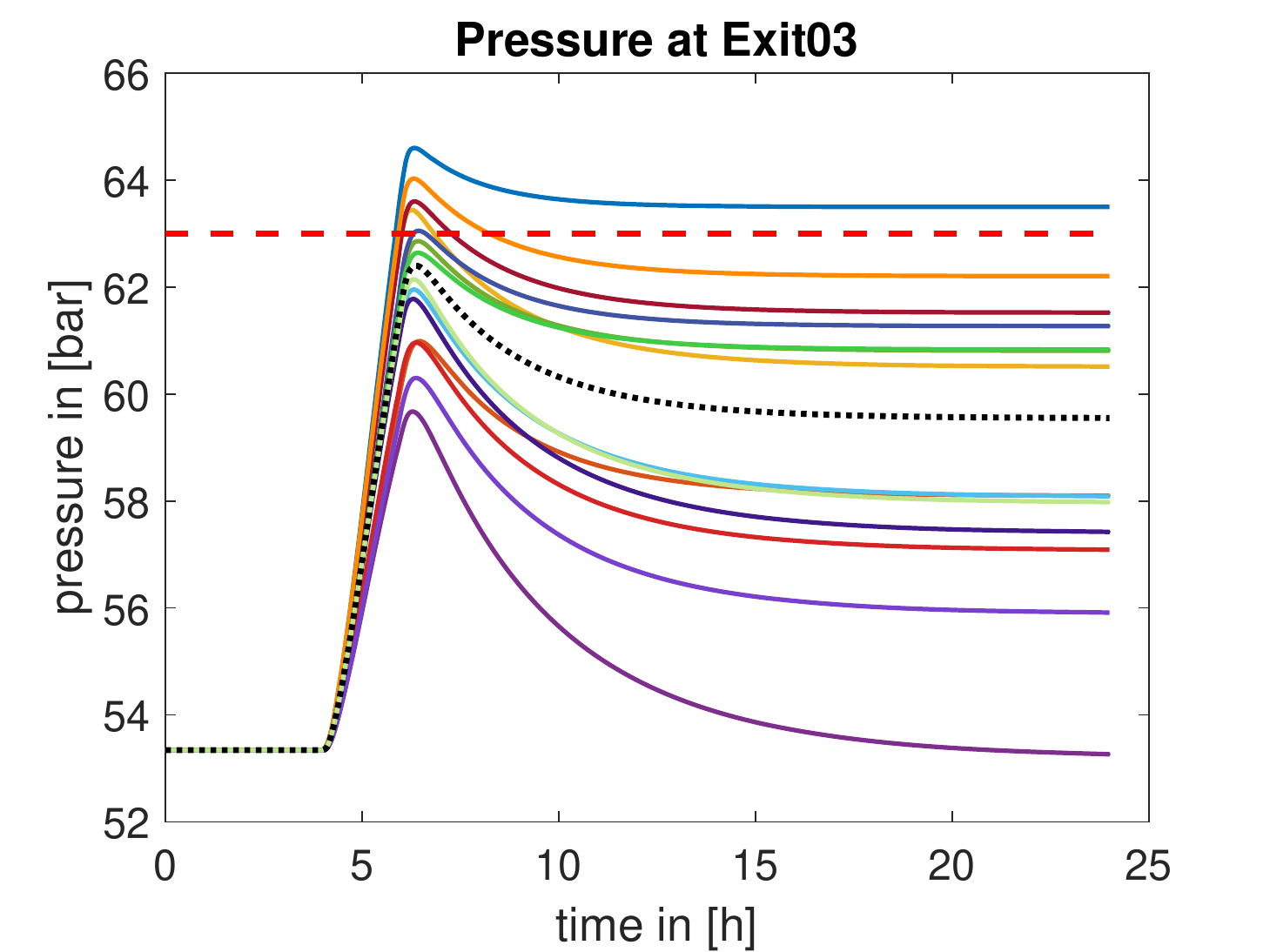}\\[2mm]
\parbox{14cm}{
\caption{\small\emph{GasLib-11}: Pressure evolution at exits E1, E2, E3, for $15$ collocation
points $y^{(j)}\in\Gamma$ chosen by the adaptive collocation method. The point
$(25.48,25.48,25.48)$ (dotted black line) corresponds to the original final state
$U_B$ with no uncertainties. The predefined pressure bounds $p^*_{min}=43\,\si{\bar}$ and
$p^*_{max}=63\,\si{\bar}$ are also plotted (red dotted lines). Obviously, these bounds
are violated by a few samples.}
\label{lds_fig:pAtExits11}
}
\end{figure}

We can also use the anisotropic Smolyak decomposition in (\ref{lds_sparse_grid_intpl})
to study the validity of the pressure bounds at the three exists $E1$, $E2$, and $E3$.
Replacing $\psi(U_h)$ by the time-dependent pressure yields
\begin{equation}
\label{lds_p11_surrogate_model}
\cI_Q[p(U_h)](y) = \sum_{\bfi\in I} \triangle^{m(\bfi)}[p(U_h)](y).\\[2mm]
\end{equation}
Exemplarily, in Fig.~\ref{lds_fig:pAtExits11}, we show the pressure curves at the
exits for $15$ collocation points
$y^{(j)}\in\Gamma$ adaptively chosen by the Smolyak algorithm for $\varepsilon=10^{-5}$.
Supposing a feasible range $[43\,\si{\bar},63\,\si{\bar}]$ for the pressure $p_{exit}$ at which the gas
should be delivered to the consumers, we are now interested in the probabilities
\begin{equation}
\PO(p_{min}<43\,\si{\bar})\quad\text{ and }\quad \PO(63\,\si{\bar}<p_{max}),
\end{equation}
with $p_{min}(y)=\min_{t\in[0,T]}p_{exit}(t,y)$ and $p_{max}(y)=\max_{t\in[0,T]}p_{exit}(t,y)$.
The surrogate model (\ref{lds_p11_surrogate_model}) allows a fast evaluation over
a sufficiently fine uniform mesh in the stochastic parameter space $\Gamma\subset\R^3$, thus
giving enough information to approximate the probability density functions of
the random variables $p_{min}$ and $p_{max}$ by a one-dimensional kernel density estimator (KDS)
\begin{equation}
KDS(x) = \frac{1}{N_sH}\sum_{i=1}^{N_s}\frac{1}{\sqrt{2\pi}}
\exp\left(-\frac{1}{2}\left(\frac{x-p(y^{(i)})}{H}\right)^2\right),
\end{equation}
for $p=p_{min},p_{max}$, where $H=1.06\,\sigma_{N_s}/N_s^{0.2}$ and $N_s=51^3$.
Observe that the bandwidth $H$ depends on the standard deviation $\sigma_{N_s}$ of
the samples as, e.g., stated and explained in \cite[Chap. 4.2]{Gramacki2018}.
\begin{figure}[!t]
\centering
\includegraphics[scale=.47]{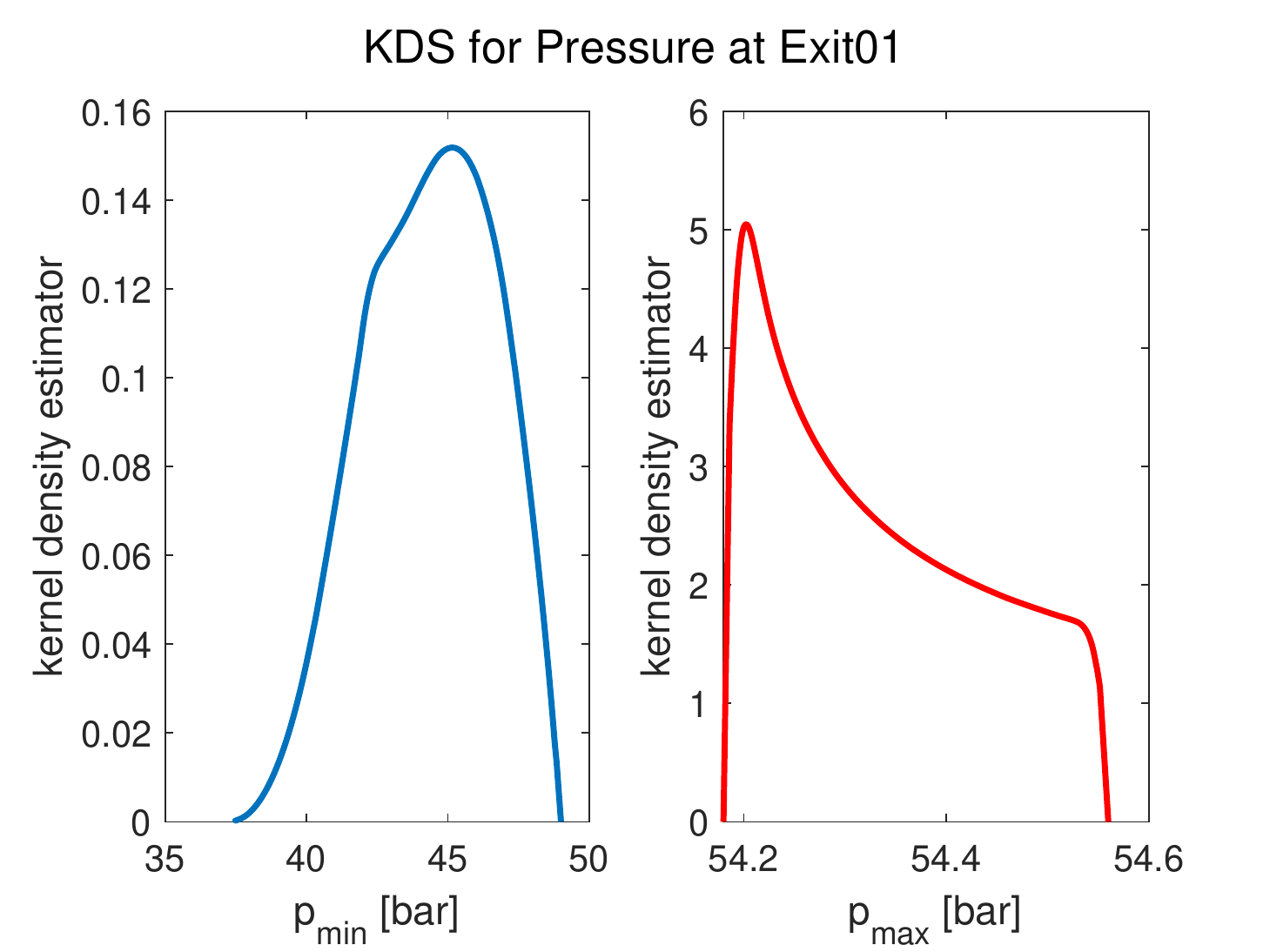}
\includegraphics[scale=.47]{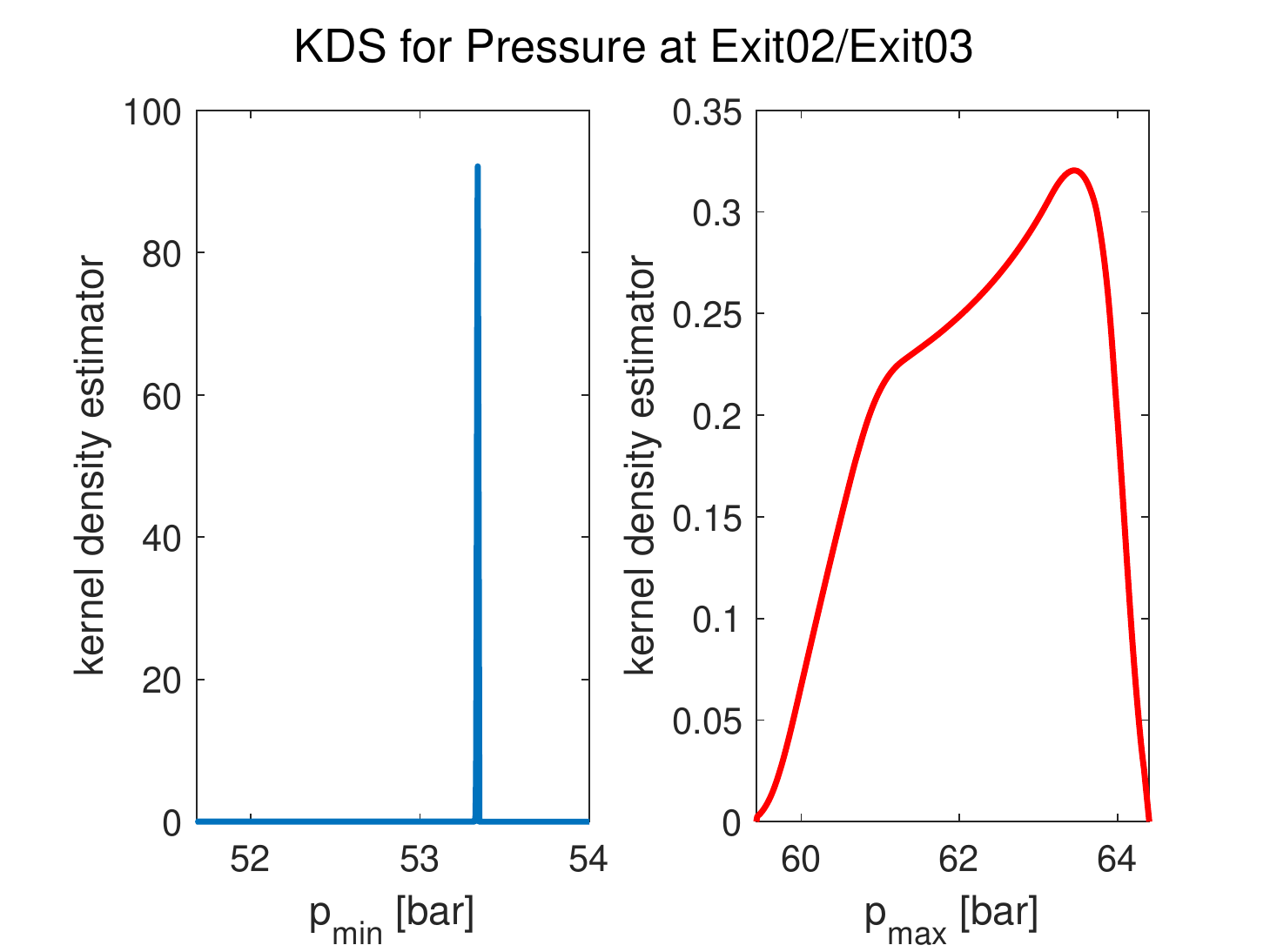}\\[2mm]
\parbox{14cm}{
\caption{\small\emph{GasLib-14}: Kernel density estimators (KDS) as
approximation of the probability density functions for the
minimum (left) and maximum (right) pressure at $E1$, $E2$ and $E3$.
Due to an inherent symmetry, the KDS for $E2$ and $E3$ are equal.}
\label{lds_fig:pdfsAtExits11}
}
\end{figure}
The corresponding KDSs are plotted in Fig.~\ref{lds_fig:pdfsAtExits11}.

From the KDSs, we calculate
\begin{equation}
\begin{array}{lcl}
\PO(p_{min}<43\,\si{\bar})&=&
\left\{
\begin{array}{ll}
0.30, & \text{ for }E1, \\
0.00, & \text{ for }E2, \\
0.00, & \text{ for }E3,
\end{array}
\right. \\[5mm]
\PO(63\,\si{\bar}<p_{max})&=&
\left\{
\begin{array}{ll}
0.00, & \text{ for }E1, \\
0.33, & \text{ for }E2, \\
0.33, & \text{ for }E3.
\end{array}
\right.
\end{array}
\end{equation}
With such information at hand, a managing operator is prepared to
react on sudden changes in the gas network with an appropriate
adaptation of the controls. It also forms the basis for probabilistic
constrained optimization, see \cite{SchusterStrauchGugatLang2020}
for more details.

\subsection{An Example with 40 Pipes}
\begin{figure}[t]
\centering
\includegraphics[scale=.7]{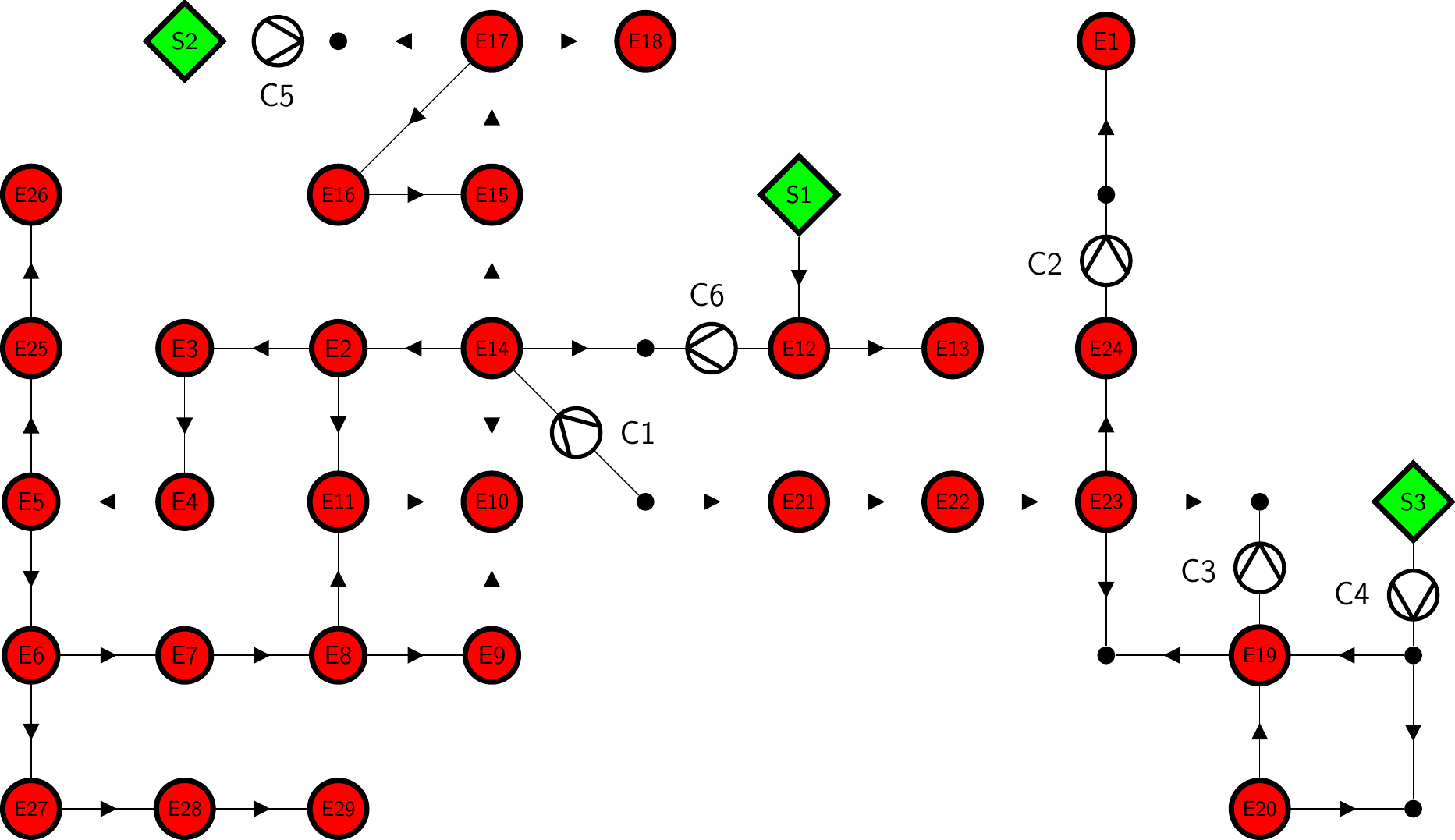}\\[5mm]
\parbox{14cm}{
\caption{\small Schematic description of the network \emph{GasLib-40} with 40 pipes,
6 compressors (C1-C6), 3 sources (green diamonds: S2, S2, S3) and 29 exits
(red circles: E1-E29). The arrows determine the
orientation of the pipes to identify the flow direction by the sign of the velocity.}
\label{lds_fig:gaslib40}
}
\end{figure}

\begin{figure}[!h]
\centering
\includegraphics[scale=0.5]{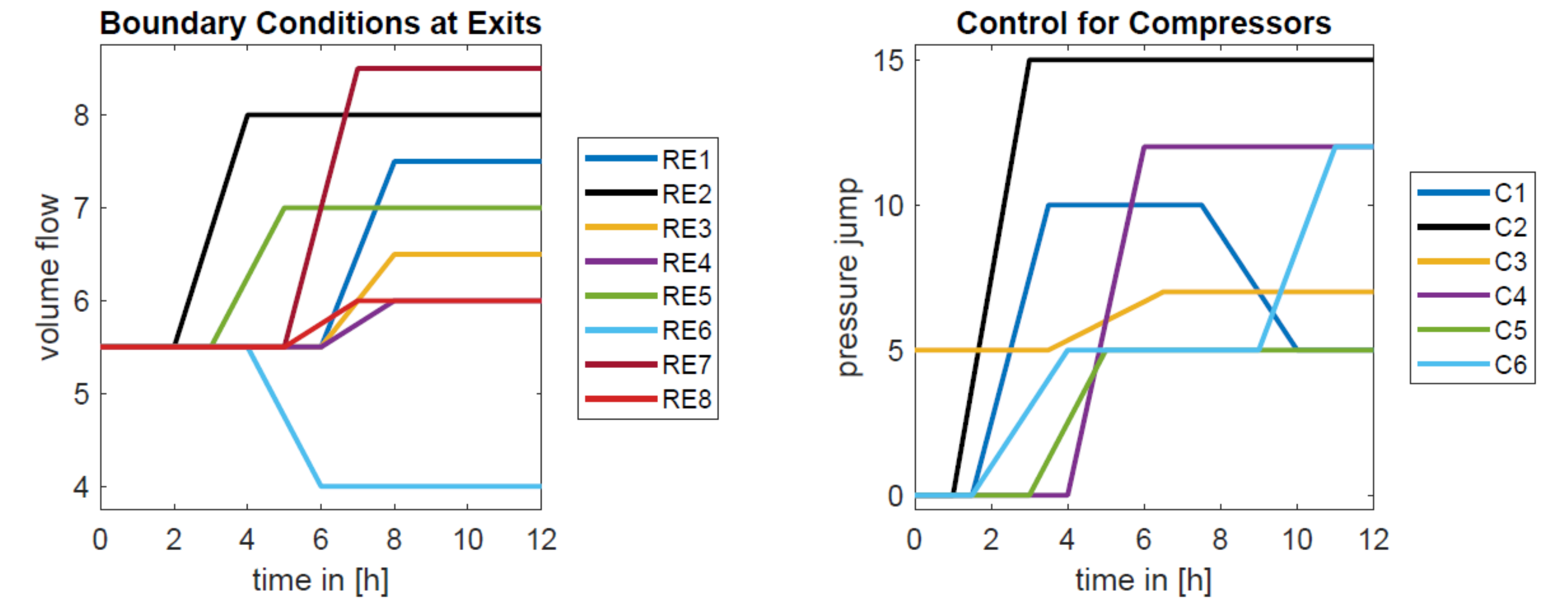}\\[2mm]
\parbox{14cm}{
\caption{\small\emph{GasLib-40}: Time-resolved boundary conditions at the exit
regions RE1-RE8 and control for the compressors C1-C6 for a smooth transition
from state $U_A$ to state $U_B$ defined in Tab.~\ref{lds_tab:4}.}
\label{lds_fig:bc_gaslib40}
}
\end{figure}

Our second example is \emph{GasLib-40}, a simplified real part of the
German Gas Network, and consists of $40$ pipes, $6$ compressor stations,
$3$ sources, and $29$ exits. Its structure is shown in Fig.~\ref{lds_fig:gaslib40}.
The exits will be clustered in $8$ different local regions (REs) with equal uptake rates
and uncertainties:
\begin{equation}
\begin{array}{l}
RE1\!=\!E1,\; RE2\!=\!E2\!-\!E11, \; RE3\!=\!E12\!-\!E13, \; RE4\!=\!E14\!-\!E18, \\
RE5\!=\!E19\!-\!E20, \; RE6\!=\!E21\!-\!E24, \; RE7\!=\!E25\!-\!E26, \; RE8\!=\!E27\!-\!E29.
\end{array}
\end{equation}
The stationary initial state $U_A$ and the final state $U_B$
are determined by the boundary conditions and controls given in Tab.~\ref{lds_tab:4}.
The temporal evolution of these values is shown in Fig.~\ref{lds_fig:bc_gaslib40}.
The computational time interval $[0h,12h]$ is split into $4$ equal subintervals.

The quantity of interest $\psi(U)$ is again defined by the specific energy consumption
of the compressors,
\begin{equation}
\psi(U(y)) = \alpha\,\sum_{c=C1,\ldots,C6} \int_{0h}^{12h} g_{c,0} + g_{c,1}\,G_c(U(y)) +
g_{c,2}\,G_c^2(U(y))\,dt
\end{equation}
with $g_{c,0}=2629$, $g_{c,1}=2.47428571429$, $g_{c,2}=1.37142857143\times 10^{-5}$
for all compressors and $G_c$ defined in (\ref{lds_gc}). The weighting factor is chosen
as $\alpha=10^{-10}$ to get values of moderate size, i.e., around $0.1$.

First, we would like to demonstrate the performance of the adaptive black box solver
ANet($\cdot,\eta_h$) for this larger network. Given the boundary conditions and controls
defined in Fig.~\ref{lds_fig:bc_gaslib40},
we always start with the initial time step $\triangle t_0=1800s$, the mesh width
$\triangle x_0=1000m$, and the simplest algebraic model
$\cM_3$. The statistics of the runs for tolerances $\eta_h=10^{-i},i=1,\ldots,5$
are summarized in Tab.~\ref{lds_tab:5}. The observed estimation process is quite reliable
and the tolerances are always satisfied. It is nicely seen that the portion of the most
detailed physical model $\cM_1$ is increasing with higher tolerances. For the last three
tolerances, we can detect $CPU\sim \eta_h^{-1}$. This was also reported for even more complex
networks in \cite{DomschkeKolbLang2020}.

Next, we model uncertainties in the exit regions RE1-RE8 by eight independent,
uniformly distributed parameters
$y=(y_1,\ldots,y_8)$, $y_i\in\cU[-1,1]$, to describe random volume flows for the state $U_B$
through
\begin{equation}
q_{REi}(y_i) = (1+0.3\cdot y_i) q_{REi}(U_B),\quad i=1,\ldots,8,
\end{equation}
where $q_{REi}(U_B)$ is the corresponding volume flow for the stationary state $U_B$ defined
in Tab.~\ref{lds_tab:4}. The parameters necessary to run the adaptive stochastic
collocation methods were determined by a few samples for low tolerances as follows:
$C_H=0.25$, $C_Y=0.1$, $s=1$, and $\mu=2$.

\begin{table}[!t]
\centering
\parbox{14cm}{
\caption{\small\emph{GasLib-40}: Boundary data for sources (S1-S3), exit regions
(RE1-RE8), and controls for
compressors (C1-C6) for initial state $U_A$ and final state $U_B$.}
\label{lds_tab:4}}
\begin{tabular}{|p{3.5cm}|cccc||cccc|} \hline
 & \multicolumn{4}{|c||}{State $U_A$} & \multicolumn{4}{|c|}{State $U_B$}\\ \hline\hline
source & S1 & S2 & S3 & & S1 & S2 & S3 & \\ \hline
pressure [\si{\bar}] & 60.0 & 53.2 & 53.2 & & 60.0 & 58.0 & 53.2 & \\ \hline\hline
exit &  RE1 & RE2 & RE3 &  RE4 & RE1 & RE2 & RE3 & RE4 \\ \hline
volume flow [\si{m^3.s^{-1}}] & 5.5 & 5.5 & 5.5 & 5.5 & 7.5 & 8.0 & 6.5 & 6.0 \\ \hline
exit &  RE5 & RE6 & RE7 & RE8 & RE5 & RE6 & RE7 & RE8 \\ \hline
volume flow [\si{m^3.s^{-1}}] & 5.5 & 5.5 & 5.5 & 5.5 & 7.0 & 4.0 & 8.5 & 6.0 \\ \hline\hline
compressor & C1 & C2 & C3 & C4 & C1 & C2 & C3 & C4 \\ \hline
pressure jump [\si{\bar}] & 0 & 0 & 5 & 0 & 5 & 15 & 7 & 12 \\ \hline
compressor & C5 & C6 & & & C5 & C6 & & \\ \hline
pressure jump [\si{\bar}] & 0 & 0 & & & 5 & 12 & & \\ \hline
\end{tabular}
\end{table}

\begin{table}[h!]
\centering
\parbox{15cm}{
\caption{\small\emph{GasLib-40}: Errors (ERR) computed from the approximate
quantity of interest
$\psi=0.1207377$ for $\eta_h=10^{-5}$, absolute value of the
sum of error estimators (EST) used in
(\ref{lds_qoi_acc}), minimum and maximum time steps $\triangle t$ and mesh
resolution $\triangle x$, distribution
of models over the pipes and computing time (CPU) for different tolerances
$\eta_h$.}
\label{lds_tab:5}}
\begin{tabular}{|c|cc|ccc|r|}
\hline
$\eta_h$ & ERR & EST &
$\begin{array}{c}\Delta t [s]\\max/min\end{array}$ &
$\begin{array}{c}\Delta x [m]\\max/min\end{array}$ &
$\begin{array}{c}\cM_1\!:\!\cM_2\!:\!\cM_3\\\text{[\%]}\end{array}$ & CPU[s] \\
\hline
$10^{-1}$ & $1.1\,10^{-2}$ & $2.2\,10^{-2}$ & 3600/1800 & 7915/767 & 0:0:100 & 3.1 \\
$10^{-2}$ & $1.3\,10^{-3}$ & $2.9\,10^{-3}$ & 3600/1800 & 7915/767 & 05:20:75 & 4.2 \\
$10^{-3}$ & $3.9\,10^{-4}$ & $1.5\,10^{-4}$ & 1800/1800 & 7915/767 & 17:51:32 & 5.3 \\
$10^{-4}$ & $1.2\,10^{-5}$ & $1.5\,10^{-5}$ & 112/112 & 7915/767 & 31:56:13 & 47.9 \\
$10^{-5}$ &           -    & $3.4\,10^{-6}$ & 28/7 & 7807/767 & 46:47:07 & 511.5 \\
\hline
\end{tabular}
\end{table}

\begin{figure}[!h]
\centering
\includegraphics[scale=.6]{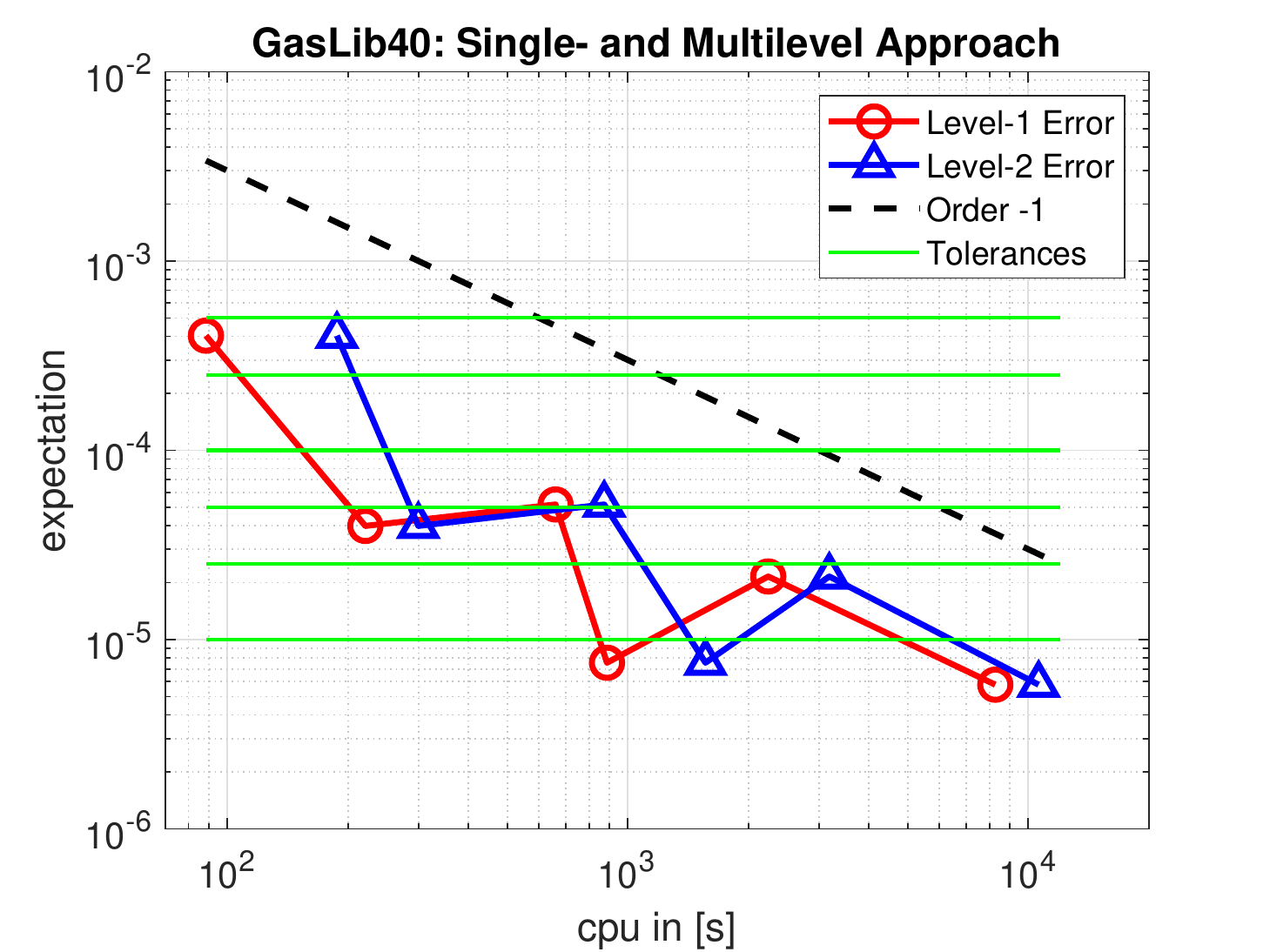}\\[2mm]
\parbox{14cm}{
\caption{\small\emph{GasLib-40}: Errors for the expected values $\EE[\psi_1^{(ML)}]$
and $\EE[\psi^{(SL)}]$ for the two-level (blue triangles), and
one-level (red circles) approach with
adaptive space-time-model discretizations for
$\varepsilon=5\times 10^{-4},2.5\times 10^{-4},10^{-4},5\times 10^{-5},2.5\times 10^{-5},10^{-5}$ (green lines). The accuracy achieved is always better than the tolerance. Both methods
deliver equal expectations since exactly the same collocation points in the
stochastic space are used.}
\label{lds_fig:slml40}
}
\end{figure}
Let us now consider a single- and two-level
approach with a reduction factor $q=0.5$ and tolerances
$\varepsilon=5\times 10^{-4},2.5\times 10^{-4},10^{-4},5\times 10^{-5},2.5\times 10^{-5},10^{-5}$.
We computed a reference solution $\EE[\psi(U)]=0.120729561141951$ with $\varepsilon=5\times10^{-7}$.
The results are shown in Fig.~\ref{lds_fig:slml40}. Both methods deliver equal values
for the expectation of $\psi(U_h)$. A closer inspection shows that only $17$ collocation points
on the finest level are sufficient to reach the desired accuracy in all runs. This also explains the observation that the two-level approach takes slightly larger computing times since the method
additionally calculates values on the coarse level. As also seen in the last example, the
numbers of samples necessary to reach the tolerances are extremely small such that the single-level
approach works already very efficient.

\section{Conclusion and Outlook}
In this study, we have applied a combination of two state-of-the-art adaptive
methods to quantify smooth
uncertainties in gas transport pipelines governed by systems of hyperbolic
balance laws of Euler type. Our in-house software tool {\sc Anaconda} and
the open-source {\sc Matlab} package \textit{Sparse Grid Kit} provide
a posteriori error estimates that can be exploited to drastically reduce the
number of degrees of freedom by using a sample-dependent strategy so that the
computational effort at each stochastic collocation point can be optimised
individually. A single-level as well as a multilevel approach have been discussed
and applied to two practical examples from the public gas library
\textit{gaslib.zib.de}.
Both strategies perform similar and quite reliable even for very high levels
of accuracy. However, we expect to see a greater potential of the multilevel
approach when facing more challenging problems in future case studies.

In contrast to Monte Carlo methods, stochastic collocation schemes provide
an access to
a global interpolant over the parameter space, which can be interpreted as
response surface approximation and used to easily calculate statistical moments and
approximate probability density functions in a postprocessing. We are
planning to incorporate
these techniques into our continuous optimization framework and thus aiming
for solving nonlinear probabilistic constrained optimization problems.

\medskip\noindent
{\bf Acknowledgements.}
The authors are supported by the Deutsche Forschungsgemeinschaft (DFG, German Research Foundation) within the collaborative research center TRR154 “Mathematical modeling, simulation and optimisation using the example of gas networks“ (Project-ID239904186, TRR154/2-2018, TP B01). We would like to thank Oliver Harbeck for making
his drawing software available to us. We have enjoyed using it to create the figures \ref{lds_fig:gaslib11} and \ref{lds_fig:gaslib40}.

\bibliographystyle{plain}
\bibliography{uqhypliterature}

\end{document}